\documentclass[a4paper,11pt,oneside,reqno]{amsart}
\usepackage{graphicx} 
\usepackage{amsfonts}
\usepackage{amsmath,amssymb}
\usepackage{amsthm}
\usepackage{xcolor}
\newtheorem{theorem}{Theorem}
\newtheorem{corollary}{Corollary}

\title[Delocalisation for $1d$ long-range  $\mathbb Z$-valued fields]{Absence of shift-invariant Gibbs states (delocalisation) 
for one-dimensional\\  $\mathbb Z$-valued fields with Long-Range interactions
}

\author{Loren Coquille}
\address{Loren Coquille, Univ. Grenoble Alpes, CNRS, Institut Fourier, F-38000 Grenoble, France}
\email{loren.coquille@univ-grenoble-alpes.fr}

\author{Aernout C.D. van Enter}
\address{Aernout C.D. van Enter, Johann Bernoulli Institute for Mathematics and Computer Science, Nijenborgh 9, 9747AG, University of Groningen, Groningen, Netherlands}
\email{a.c.d.van.enter@rug.nl}

\author{Arnaud Le Ny}
\address{Arnaud Le Ny, LAMA UMR CNRS 8050, UPEC, Universit\'e Paris-Est, 61 Avenue du G\'en\'eral de Gaulle,  94010 Cr\'eteil cedex, France}
\email{arnaud.le-ny@u-pec.fr}

\author{Wioletta M. Ruszel}
\address{Wioletta M. Ruszel, Utrecht University, Hans Freudenthalgebouw
Budapestlaan 6, 3584 CD Utrecht, The Netherlands}
\email{w.m.ruszel@uu.nl}

\date{June 2024}

\begin{document}
\thanks{\textit{Acknowledgments.} We thank Christophe Garban for a helpful feedback. \\
Research of A.L.N. and L.C. have been supported by the CNRS IRP (International Research Project) EURANDOM “Random
Graph, Statistical Mechanics and Networks” and by the LabEx PERSYVAL-
Lab (ANR-11-LABX-0025-01) funded by the French program Investissement
d’avenir. W.M.R. is funded by Vidi grant VI.Vidi.213.112
from the Dutch Research Council.
}
\maketitle
\begin{abstract}
We show that a modification of the proof of our paper~\cite{CvELNR}, in the spirit of \cite{FP81}, shows delocalisation in the long-range Discrete Gaussian Chain, and generalisations thereof, for any decay power $\alpha>2$ and at all temperatures.
The argument proceeds by contradiction: any shift-invariant and localised measure (in the $L^1$ sense), is a convex combination of ergodic localised measures. But the latter cannot exist: on one hand, by the ergodic theorem, the average of the field over growing boxes would be almost surely bounded ; on the other hand the measure would be absolutely continuous with respect to its height-shifted translates, as a simple relative entropy computation shows. This leads to a contradiction and answers, in a non-quantitative way, an open question stated in a recent paper of C.~Garban  \cite{G}.
\end{abstract}

\footnotesize




\normalsize

\section{Introduction}

The long-range Discrete Gaussian Chain (DGC) is a $\mathbb{Z}$-valued version of the $1d$ long-range Ising model, where the height variables $\phi_i \in \mathbb{Z}$ at site $i \in \mathbb{Z}$ form a height configuration $\phi$ in the state-space $\mathbb{Z}^\mathbb{Z}$, which is sampled according to the Gibbs measure with formal Hamiltonian

\begin{equation} \label{DGC}
	H(\phi) :=\sum_{i \neq j \in\mathbb Z} J_\alpha(|i-j|) \cdot |  \phi_i - \phi_j |^2
\end{equation}

with long-range coupling constants $J_\alpha(k) \sim k^{-\alpha}$,
the most classical choice being $J_\alpha(k)=  k^{-\alpha}$. It has been introduced as an effective model for an interface between two ordered phases in $2d$ long-range Ising models with the
same polynomially decaying coupling constants, and was studied in \cite{KH82,FZ91, Sheffield-Ast,Velenik06, CvELNR}.

The localisation/delocalisation of models with Hamiltonians (\ref{DGC}) is analogous to the same question for $2d$  nearest neighbour integer-valued GFF (also called
Discrete Gaussian model) with $\phi : \mathbb{Z}^2 \longrightarrow \mathbb{Z}$ whose delocalisation at high temperature has been proved in \cite{FS81,Lam22}, see also \cite{LO}. \\

The behaviour of the long-range DGC with Hamiltonian \eqref{DGC} is as follows:
\begin{itemize}
	\item For $\alpha\in(1,2)$ the field $\phi$ is localised at any inverse temperature $\beta$, as recently proved in \cite{G}.
	\item For $\alpha=2$, there is a localisation/delocalisation phase transition, which was originally proved by Fröhlich and Zegarlinski \cite{FZ91}. More precisely, the localisation at low temperature was proved in \cite{FZ91} and the delocalisation at high temperature was proved in \cite{KH82,FZ91}.
	\item For $\alpha>2$, Fröhlich and Zegarlinski \cite{FZ91} conjecture that fluctuations are of order $n^{\min\{\alpha-2,1\}/2}$ at all temperatures, and give a heuristical energy-entropy argument.
	\item For $\alpha\in[2,3)$, the paper \cite{G} proposes a new proof of delocalisation at high temperature, relying on \cite{FS81}, which goes much beyond the above conjecture and exhibits the phenomenon of "invisibility of the integers": an invariance principle towards a $H$-fractional Brownian motion of Hurst index $\frac{\alpha-2}2$ holds, and this limiting Gaussian field is the same as the one arising from the \emph{real}-valued Gaussian Chain. The diffusivity constant is also computed.
	\item For $\alpha>3$, there is delocalisation at all temperatures and fluctuations follow a standard central limit scaling of order $n^{1/2}$, see \cite{G}.
\end{itemize}
\vspace{.2cm}
The following open problem is stated in \cite{G} as {\bf Open Problem 5}:

\textit{Prove the delocalisation of the discrete Gaussian chain at all temperatures when $\alpha \in (2,3)$, possibly with quantitative estimates on the variance. A heuristics is discussed in \cite{FZ91}, but as far as we know, a proof is still missing.}\\

In the work \cite{CvELNR}, the four present authors studied the $2d$ long-range Ising model, and proved absence of Dobrushin states (extremal non-transation invariant Gibbs measures arising from Dobrushin boundary conditions) for any $\alpha>2$ and at any temperatures, by means of Araki's relative entropy method from~\cite{Araki,FP81,FP86}.\\

The goal of the present note is to show that the same method leads to a (non-quantitative) proof of delocalisation of the long-range DGC or SOS model, and generalisations thereof, at all temperatures for any $\alpha>2$. Quantitative estimates on the delocalisation should be derivable, providing some bounds on the fluctuation exponent, which is conjectured in \cite{FZ91} to be equal
to $\min\{\alpha-2, 1\}/2$ as stated above, but in this short note we restrict ourselves to  the following rigorous qualitative statement. \\

{Note that such a simple proof would (and should) not work in the two-dimensional case, because the relative entropy between the field and its vertical translate in a box of size $n$ by $n$ will grow as fast as the perimeter of the box\footnote{The constant $C_3$ in Equation \eqref{cst} of the proof below would be replaced by $C_3n$.}, which is not bounded in dimension 2. 
	This is consistent with the fact that at low temperature,  the $2d$ nearest neighbor SOS model localises, so the long-range models should be localised also (even more by inequalities). Delocalisation at high temperature when $\alpha >4$ is the plausible conjecture (see \cite{G}, Open problem 7) but we expect the proof to be much harder.}


\section{Non-quantitative proof of delocalisation}

Consider $\mathbb Z$-valued fields on $\mathbb{Z}$ under potential $V(x)=|x|^p$ with $p\in[1,2]$ and long-range polynomially decaying interactions $J_\alpha(x)\sim|x|^{-\alpha}$.
The Hamiltonian corresponding to this class of models is : 
\begin{equation}\label{DGC-SOS}
	H(\phi)= - \sum_{i\neq j\in\mathbb Z} J_{ij} V(\phi_i - \phi_j),
\end{equation}
where $J_{ij} = J_\alpha(|i-j|)$ (positivity is not assumed)
and $\phi_i$ takes values in $\mathbb Z$.
The model makes sense for $\alpha>1$.
The DGC corresponds to $V(x)=x^2$, and the SOS model corresponds to $V(x)=|x|$. For general $p$, the model is also called $p$-SOS.


Let $\mu^\omega_\Lambda$ denote the finite-volume measure with Hamiltonian \eqref{DGC-SOS} with {height-$\omega$} boundary condition outside the box $\Lambda\Subset\mathbb Z$. More precisely, for $\beta>0$,
$$\mu^\omega_\Lambda(\phi)\propto\exp\left(-\beta H(\phi)\right)1_{\phi_i\in\mathbb Z\;\forall i\text{ and }\phi_i=\omega\;\forall i\in\Lambda^c}.$$
In the case $\omega\equiv0$ we write $\mu^0_\Lambda$.
A Gibbs measure for the Hamiltonian \eqref{DGC-SOS} is a probability measure $\mu$ on $\Omega=\mathbb Z^{\mathbb Z}$ such that for any $\Lambda\Subset\mathbb Z$, we have
\begin{equation}\label{DLR}
	\mu(\phi)=\int d\mu(\omega)\mu^\omega_\Lambda(\phi).
\end{equation}

Expectations under a measure $\nu$ are written $\langle\cdot\rangle_\nu$.
We say that a measure $\nu$ is shift-invariant if $\nu(\phi)=\nu(\tau(\phi))$ for any translation $\tau:z\mapsto z+v$, $v\in\mathbb Z$.
The set of shift-invariant Gibbs measures is denoted by $\mathcal G_{si}$.\\
We write $\nu_t(\phi)=\nu(\phi+t)$, when we shift the heights by an amount $t$.
\begin{theorem} \label{THM}
	For any $\beta>0$ and $\alpha>2$, delocalisation holds for the models with Hamiltonian \eqref{DGC-SOS} in the following sense:
	\begin{enumerate}
		\item Let $p\in\{1,2\}$, i.e. consider either the SOS model or the DGC. \\
		There does not exist any shift-invariant Gibbs measure $\mu\in\mathcal G_{si}$ such that {localisation holds in the $L^1$ sense}
		\begin{equation}\label{LOC-1}
			{\langle|\phi_0|\rangle}_\mu < {\infty}.
		\end{equation}
		\item Let $p\in(1,2)$. 
		There does not exist any shift-invariant Gibbs measure $\mu\in\mathcal G_{si}$ such that {localisation holds in the $L^2$ sense}
		\begin{equation}\label{LOC-2}
			{\langle\phi_0^2\rangle}_\mu < {\infty}.
		\end{equation}
	\end{enumerate}
\end{theorem}
{
	\begin{corollary}\label{COR}
		For any $\beta>0$, $\alpha>2$, also finite-volume delocalisation holds for the ferromagnetic models (i.e. $J_{ij}\geq0$) with Hamiltonian \eqref{DGC-SOS}, in the following sense :\\
		For any sequence of boxes $(\Lambda_n)_n\uparrow\mathbb Z$ it holds 
		\begin{equation}\label{finite-volume-deloc}
			{\langle\phi_0^2\rangle}^0_{\Lambda_n}\to\infty\quad
			\text{ as }n\to\infty.
		\end{equation}
	\end{corollary}
	\begin{proof}[Proof of Corollary \ref{COR}] 
		Lammers and Ott proved the FKG inequality for the absolute value of the field (see \cite{LO}, Theorem 2.8) 
		which implies, in the context of \emph{planar} interactions,  the equivalence of the two notions of delocalisation, namely (D1) delocalisation in finite volume in the sense of \eqref{finite-volume-deloc}, and (D2) absence of shift-invariant Gibbs measures (see \cite{LO}, Theorem 2.7).\\
		We notice that their proof of the implication (D2)$\Rightarrow$(D1) holds without the planarity assumption on the interactions.\\
		Indeed, following Section 7 of \cite{LO}, one notices that $\phi_i=0\Leftrightarrow|\phi_i|=0$. The FKG lattice condition for the absolute value of the field (for which zero is minimal) and the DLR property \eqref{DLR} imply that for any $\Lambda\subset\Lambda'\subseteq\mathbb Z$, the stochastic domination $|\phi|\preceq|\phi'|$ holds under $\mu^0_\Lambda$ and $\mu^0_{\Lambda'}$. In particular, the variances $\langle\phi_i^2\rangle^0_{\Lambda_n}=\langle|\phi_i|^2\rangle^0_{\Lambda_n}$ are increasing in $n$, for any $i\in\mathbb Z$, and the limiting variances $$\lim_{n\to\infty}{\langle\phi_i^2\rangle}^0_{\Lambda_n}$$ exist in $[0,\infty]$ and  coincide for all $ i\in\mathbb Z$. \\
		Assume the limiting variance is finite. We can thus follow \cite{LO} to construct a shift-invariant localised Gibbs measure, which cannot exist by Theorem \ref{THM}, and leads to a contradiction, proving Corollary \ref{COR}.
	\end{proof}
}
{Note that we do not} exclude any (subsequential) weak limit of the $\mu^0_{\Lambda_n}$ not to be extremal, nor exclude it to be a mixture of exclusively non-shift invariant extremal measures.
Such exclusions can be proved in the case of planar interactions, see the thesis of Sheffield \cite{Sheffield-Ast}.

\begin{proof}[Proof of Theorem \ref{THM}]The argument is based on \cite{FP81,CvELNR}, and proceeds by contradiction: if we suppose $\mu$ is shift-invariant and localised (\eqref{LOC-1} holds), then it is a convex combination of ergodic localised measures. But an ergodic localised measure $\nu$ cannot exist: on one hand, the average of the field over growing boxes would be $\nu$ almost surely bounded by ergodicity and by our under $L^1$ assumption \eqref{LOC-1}; on the other hand $\nu$ would be absolutely continuous with respect to its vertical translate $\nu_t$ as a simple relative entropy computation shows. This leads to a contradiction for $t$ large enough. \\
	\noindent More precisely, fix $\beta>0$ and suppose there exists $\mu\in\mathcal G_{si}$ which is localised:
	\begin{equation}\label{LOC}
		\exists C\in[0,\infty)\text{ such that }\langle|\phi_0|\rangle_\mu < C.
	\end{equation}
	Then $\mu$ has a unique decomposition $w_\mu$ onto shift-invariant extremal Gibbs measures (see \cite{Georgii}, Theorem 14.17):
	\begin{equation*}
		\mu=\int_{ex\mathcal G_{si}}\nu\;w_\mu(d\nu).
	\end{equation*}
	Moreover, this decomposition is supported on localised measures. Indeed, as
	$
	\langle|\phi_0|\rangle_\mu =
	\int_{ex\mathcal G_{si}}{\langle|\phi_0|\rangle}_\nu w_\mu(d\nu)< C,
	$
	then for $w_\mu$ almost every $\nu\in ex\mathcal G_{si}$, there exists some $C_\nu<\infty$ such that
	\begin{align}\label{LOCnu}
		{\langle|\phi_0|\rangle}_\nu\leq C_\nu.
	\end{align}
	We now prove that an extremal shift-invariant localised Gibbs measure $\nu$ for the Hamiltonian \eqref{DGC-SOS} cannot exist if $\alpha>2$.\\
	Assume $\nu\in ex\mathcal G_{si}$ is such that \eqref{LOCnu} holds.

	On the one hand, any $\nu\in ex\mathcal G_{si}$ is ergodic (see \cite{Georgii} Theorem 14.15) and thus the field is strongly localised under $\nu$, in the following sense. \\
	Let $M_n=\frac1{2n}\sum_{|i|\leq n}\phi_i$. By the ergodic theorem, under hypothesis \eqref{LOC-1}, 
	\begin{align}\label{strong-LOC}
		\nu\big(M_n\to {\langle\phi_0\rangle}_\nu\big)=1,
	\end{align}
	and thus for any $t\in\mathbb Z$,
	\begin{align}\label{strong-LOC-t}
		\nu_t\big(M_n\to {\langle\phi_0\rangle}_\nu+t\big)=1.
	\end{align}
	{In particular, for any $t\neq0$ the measures $\nu$ and $\nu_t$ are singular with respect to each other: 
		\begin{equation}\label{sing}
			\forall t\neq 0, \quad\nu\perp\nu_t.
	\end{equation}}
	
	On the other hand, the relative entropy between $\nu$ and $\nu_t$ is bounded if $\alpha>2$. To prove this, we follow Section 4 of \cite{FP81} with a simplified local transformation of the measure $\nu$.
	Recall that the relative entropy between two measures $\mu$ and $\nu$, denoted $RE(\mu | \nu)$, is defined to be
	$$
	RE(\mu | \nu):= {\langle \log( {d \mu}/{d \nu} )\rangle}_{\mu}
	$$
	when $ \mu \ll \nu$, and $+\infty$ otherwise.
	Let $t\in\mathbb N$ be fixed, and define "the step of height $t$" as
	\begin{align*}
		a_i^n(t)&=\left\{
		\begin{array}{ll}
			t &\text{if }|i|<n\\
			0&\text{if }|i|\geq n.
		\end{array}
		\right.
	\end{align*}
	Define the transformation $T_{t,n}$ ("adding the step of height $t$ in a box of size $n$") on $\Omega=\mathbb Z^{\mathbb Z}$ by setting
	$$(T_{t,n}\phi)_i=\phi_i+a_i^n(t).$$
	Let $\nu_{t,n}$ be the transformed Gibbs measure defined by
	$${\langle f\rangle}_{\nu_{t,n}}=\int\nu(d\phi)f(T_{t,n}^{-1}\phi).$$
	Then $\nu_{t,n}$ is absolutely continuous with respect to $\nu$, with Radon-Nykodim derivative given by\footnote{this follows from the fact that $\nu$ is a Gibbs measure.}
	$$\frac{d\nu_{t,n}}{d\nu}(\phi)=\exp(-\beta\sum_{i\neq j\in\mathbb Z} J_{ij}(V(T_{t,n}\phi_i-T_{t,n}\phi_j)-V(\phi_i-\phi_j))).$$
	To compute the relative entropy, note that
	\begin{align}\label{1term}
		&{\langle V(T_{t,n}\phi_i-T_{t,n}\phi_j)-V(\phi_i-\phi_j)\rangle}_\nu\nonumber\\
		&={\langle|\phi_i+a_i^n(t)-\phi_j-a_j^n(t)|^p-|\phi_i-\phi_j|^p\rangle}_\nu.
	\end{align}
	In the case $p=2$, we get
	\begin{equation*}
		\eqref{1term}=2(a_i^n(t)-a_j^n(t)){\langle\phi_i-\phi_j\rangle}_\nu+(a_i^n(t)-a_j^n(t))^2=(a_i^n(t)-a_j^n(t))^2.
	\end{equation*}
	In the cases $p\in[1,2)$, we have the upper bound\footnote{using $|a+b|^p\leq 2^{p-1}(|a|^p+|b|^p)$ for any $p\geq1$.} 
	\begin{align*}
		\eqref{1term}\leq 
		(2^{p-1}-1){\langle|\phi_i-\phi_j|^p\rangle}_\nu
		+2^{p-1}|a_i^n(t)-a_j^n(t)|^p.
	\end{align*}
	The first term vanishes in the case $p=1$, and is bounded in the remaining cases $p\in(1,2)$ under hypothesis \eqref{LOC-2}, as $L^2(\nu)\subset L^p(\nu)$. \\
	Altogether, assuming \eqref{LOC-1} for $p\in\{1,2\}$ and \eqref{LOC-2} for $p\in(1,2)$, we have
	\begin{align}\label{cst}
		RE(\nu|\nu_{t,n})&=
		\big\langle\sum_{i\neq j\in\mathbb Z} J_{ij}(V(T_{t,n}\phi_i-T_{t,n}\phi_j)-V(\phi_i-\phi_j)){\big\rangle}_\nu\nonumber\\
		&\leq
		\sum_{|i|\leq n}\sum_{j\in\mathbb Z}|J_{ij}|(C_1+C_2(a_i^n(t)-a_j^n(t))^p)\nonumber\\
		&=
		\sum_{|i|\leq n}\sum_{|j|>n} |J_{ij}|(C_1+C_2t^p)
		\leq (C_3 + C_4 n^{2-\alpha})(C_1+C_2t^p)\\
		&<C_5+C_6t^p \text{ if } \alpha>2.\nonumber
	\end{align}
	where the above constants are positive and finite, they may depend on $p$ but do not depend on $n$.

	By the lower-semicontinuity property of the relative entropy\footnote{$RE(\mu|\nu)\leq\lim_{n\to\infty}RE(\mu_n|\nu_n)$ whenever $\mu_n\to\mu$ and $\nu_n\to\nu$ weakly.}, we get the same upper bound as above on the relative entropy between the infinite-volume Gibbs measures $\nu$ and $\nu_t$ (since $\nu_t$ is a weak limit of the $\nu_{t,n}$):
	\begin{align*}
		RE(\nu | \nu_t)\leq \lim_{n\to\infty}RE(\nu | \nu_{t,n})
		&<C_5+C_6t^p \text{ if } \alpha>2.
	\end{align*}
	
	We deduce  that, for any $t\in\mathbb Z$, $\nu$ and $\nu_t$ are absolutely continuous with respect to each other if $\alpha>2$.
	{This is in contradiction with \eqref{sing} and proves the Theorem.}
\end{proof}


\end{document}